\documentclass[twocolumn]{autart}    %

\usepackage[utf8]{inputenc}
\usepackage{siunitx}
\usepackage{amsmath}
\usepackage{amsfonts}
\usepackage{graphicx}
\usepackage{float}
\usepackage{bm}
\usepackage{color,xcolor}
\usepackage{mathtools}
\mathtoolsset{showonlyrefs}

\usepackage{natbib}

\usepackage{natbib}

 \usepackage{hyperref}
 \hypersetup{
 	hidelinks,
     colorlinks=false,
     linkcolor=blue,
     filecolor=blue,      
     urlcolor=blue,
 }

\newtheorem{theorem}{Theorem}%
\newtheorem{lemma}[theorem]{Lemma}

\newtheorem{proposition}{Proposition}
\newtheorem{problem}{Problem}
\newtheorem{assumption}{Assumption}
\theoremstyle{definition}
\newtheorem{definition}[theorem]{Definition}
\newtheorem{example}{Example}
\newtheorem{remark}{Remark}

\newcommand{\sign}[1]{\mbox{sign}(#1)}
\newcommand{\sgn}[1]{\lfloor#1\rceil}
\newcommand{\real}[1]{\mbox{Re}(#1)}
\renewcommand{\exp}[1]{\mbox{exp}(#1)}
\newcommand{\diag}[1]{\mbox{diag}(#1)}

\begin{document}

\begin{frontmatter}

\title{
An arbitrary-order exact differentiator with predefined convergence time bound for signals with exponential growth bound\thanksref{footnoteinfo}
} %

\thanks[footnoteinfo]{This paper was not presented at any IFAC 
meeting. 
Corresponding author: R.~Aldana-L\'opez.
}

\thanks[footnoteinfo2]{\textcolor{red}{This is the accepted manuscript version for: Gómez-Gutiérrez D., Aldana-López R., Seeber R., Angulo M.T., Fridman L. ``An arbitrary-order exact differentiator with predefined convergence time bound for signals with exponential growth bound". Automatica. 2023; DOI: 10.1016/j.automatica.2023.110995. 
\textbf{Please cite the publisher's version}. For the publisher's version and full citation details, see:
\url{https://doi.org/10.1016/j.automatica.2023.110995}.
© 2023. This manuscript version is made available under the CC-BY-NC-ND 4.0 license \url{https://creativecommons.org/licenses/by-nc-nd/4.0/}}}

\author[Intel,TECMM]{David~G\'omez-Guti\'errez}\ead{david.gomez.g@ieee.org},
\author[Zaragoza]{Rodrigo~Aldana-L\'opez\thanksref{footnoteinfo2}}\ead{rodrigo.aldana.lopez@gmail.com}, 
\author[Austria]{Richard Seeber}\ead{richard.seeber@tugraz.at}, 
\author[UNAM]{Marco Tulio Angulo}\ead{mangulo@im.unam.mx }, 
\author[UNAMMX]{Leonid Fridman}\ead{lfridman@unam.mx}

\address[Intel]{Intel Tecnolog\'ia de M\'exico, Intel Labs, Intelligent Systems Research Lab, Jalisco, Mexico.}

\address[TECMM]{Tecnológico Nacional de México, Instituto Tecnológico José Mario Molina Pasquel y Henríquez, Unidad Académica Zapopan, Jalisco, Mexico.}

\address[Zaragoza]{University of Zaragoza, Departamento de Informatica e Ingenieria de Sistemas (DIIS), Zaragoza, Spain.}

\address[Austria]{Graz University of Technology, Institute of Automation and Control, Christian Doppler Laboratory for Model Based Control of Complex Test Bed Systems, Graz, Austria. }

\address[UNAM]{CONACyT - Institute of Mathematics, Universidad Nacional Autonoma de México,  Juriquilla, Mexico.}

\address[UNAMMX]{Universidad Nacional Autonoma de México, Engineering Faculty, Ciudad Universitaria, Mexico}
          
\begin{keyword}                           %
Fixed-time stability, predefined-time, prescribed-time, unknown input observers, online differentiators              %
\end{keyword}                             %

\begin{abstract}                          %
There is a growing interest in differentiation algorithms that converge in fixed time with a predefined Upper Bound on the Settling Time (UBST). However, existing differentiation algorithms are limited to signals having an $n$-th order Lipschitz derivative. Here, we introduce a general methodology based on time-varying gains to circumvent this limitation, allowing us to design $n$-th order differentiators with a predefined UBST for the broader class of signals whose $(n+1)$-th derivative is bounded by a function with bounded logarithmic derivative. Unlike existing methods whose time-varying gain tends to infinity, our approach yields a time-varying gain that remains bounded at convergence time. We show how this last property maintains exact convergence using bounded gains when considering a compact set of initial conditions and improves the algorithm's performance to measurement noise. 
\end{abstract}

\end{frontmatter}

\section{Introduction}

The design of arbitrary order exact differentiators is instrumental in solving a wide range of estimation and control problems \citep{Levant1998RobustTechnique,Levant2003,Levant2019,Sanchez2016construction,Reichhartinger2018AnGains,Reichhartinger2017AMatlab/Simulink, Fridman2008Higher-orderSystems,Fridman2011High-orderInputs,Rios2018ASystems,Rios2015FaultApproach,Alwi2011FaultModes,FerreiraDeLoza2015OutputIdentification,Shtessel2014ObservationObservers,Imine2011ObservationHOSM-observers}. For problems with time constraints,  it is very useful if the user can a-priori select an Upper Bound on the Settling-Time (UBST) of the differentiator. For example, in output-based control problems, we can use differentiators with small enough UBST to build an observer that converges before the system trajectories exit a safe compact set. After the differentiator convergence, we can turn on a controller designed assuming exact knowledge of the entire system state.

For the classic exact differentiators by~\cite{Levant2003} and \cite{Levant2019}, selecting a-priori the UBST requires knowing the set of initial conditions, which is unfeasible for most applications. The recently developed notion of fixed-time stability can circumvent this fundamental problem by requiring a uniform UBST for all initial conditions 
\citep{Cruz2018}. However, fixed-time stability has been typically characterized using homogeneity properties~\citep{Andrieu2008}, without explicitly computing a UBST (see, e.g., ~\citep{Angulo2013RobustDifferentiator}). \cite{Cruz-Zavala2011} and \cite{Seeber2021robust} have proposed first-order differentiators with explicit UBST, but these bounds are very conservative and lead to differentiation errors that are larger than necessary. Additionally, the above methods assume that the $n-$th derivative of the noise-free input is Lipschitz. This assumption is reasonable when the system is known to converge toward a bounded invariant set, such as chaotic systems with state-dependent disturbances~\citep{Gomez-Gutierrez2017OnModulation}. However, this assumption is not satisfied if the system is unstable, even if it has linear dynamics. Importantly,  note that the trajectories of unstable linear or Lipschitz nonlinear systems cannot grow faster than an exponential ~\citep{Bejarano2011Finite-timeInputs,Rodrigues2018GlobalFeedback}. Therefore, for applications,  it is reasonable to assume signals with a bounded logarithmic derivative~\citep{Oliveira2017GlobalControl,Rodrigues2018GlobalFeedback}.

Previous works have proposed exact finite-time convergent differentiators for signals whose $(n+1)-$th derivative has a time-varying bound with bounded logarithmic derivative~\citep{Levant2012ExactDerivatives,Levant2018GloballyGains,Moreno2018ExactGains}. However, these differentiators do not have fixed-time convergence. Moreover, even for bounded initial conditions, there is no methodology to select the desired UBST for differentiators in this class of signals. 

Recently, \cite{Holloway2019} proposed an alternative design methodology of differentiators for polynomial signals of $n-$th order using Time-Varying Gains (TVG). This methodology produces an algorithm with prescribed-time convergence ---that is, for every nonzero initial condition, the differentiator converges precisely at the time the user prescribes. Prescribing the convergence time is a significant advantage 
compared to differentiators having a conservative estimation of their UBST. However, this method requires that the TVG tends to infinity at the convergence time, making its application challenging under measurement noise or limited numerical precision. Although several workarounds have been suggested to circumvent this challenge and maintain a bounded TVG,  they no longer obtain a zero differentiation error. Furthermore, with such workarounds, the error at the prescribed time grows linearly with the initial condition

To fill the above gaps, here we ``redesign" the classic differentiator of ~\cite{Levant2018GloballyGains} by using TVGs to obtain a new exact differentiator for signals with bounded logarithmic derivative, where the desired UBST is set a priori as one parameter of the algorithm. Unlike~\cite{Menard2017, Cruz-Zavala2011, Seeber2021robust}, our approach yields an arbitrarily tight UBST. Furthermore, unlike~\cite{Holloway2019}, which we consider the closest approach, we obtain exact convergence with a bounded TVG at the settling-time instant. This last property is important because it allows us to introduce workarounds that maintain exact convergence with a bounded TVG when considering a compact set of initial conditions and maintain accuracy under measurement noise comparable to the accuracy of the original differentiator. Compared to our previous result~\citep{Aldana-Lopez2021AGains} that considers only a first-order differentiator,  here we consider an arbitrary-order differentiator and extend the class of TVG we use. We also show numerically that these features allow better transient behavior. Moreover, we present the accuracy analysis under measurement noise under the proposed workaround. 

The rest of the manuscript is organized as follows. Section~\ref{Sec:Prelim} introduces the problem statement and some preliminaries. In Section~\ref{Sec:Main}, we present the main result. In Section~\ref{sec:practical}, we introduce a workaround to maintain a bounded TVG and an analysis of its accuracy under measurement noise. Section~\ref{Sec:Discussion} discusses the main features of our redesign methodology, contrasting it to other state-of-the-art algorithms. Finally, in Section~\ref{Sec:Conclu}, we present conclusions and future work. Proofs are collected in the Appendix. 

\noindent \textbf{Notation:}
$\Bar{\mathbb{R}}=\mathbb{R}\cup\{-\infty,\infty\}$, $\mathbb{R}_+=\{x\in\mathbb{R}\,:\,x\geq0\}$ and $\Bar{\mathbb{R}}_+=\mathbb{R}_+\cup\{\infty\}$.
For a signal $y: \mathbb R_+ \rightarrow \mathbb R$, we denote by $y^{(i)}$ its $i$-th derivative. 
For $x\in\mathbb{R}$, $\sgn{x}^\alpha = |x|^\alpha \mbox{sign}(x)$, if $\alpha\neq0$ and $\sgn{x}^\alpha = \mbox{sign}(x)$ if $\alpha=0$. 
For a function $\phi:\mathcal{I}\to\mathcal{J}$, its reciprocal $\phi(\tau)^{-1}$, $\tau\in\mathcal{I}$,  is such that $\phi(\tau)^{-1}\phi(\tau)=1$ and its inverse function $\phi^{-1}(t)$, $t\in\mathcal{J}$, is such that $\phi(\phi^{-1}(t))=t$. Given a matrix $A$, $A^T$ is its transpose. For $v\in\mathbb{R}^{n\times 1}$, $\|v\|=\sqrt{v
^Tv}$.
$\mathcal{U}:=[u_{ij}]\in\mathbb{R}^{(n+1)\times(n+1)}$, where $u_{ij}=1$ if $j=i+1$,  $u_{ij}=0$, otherwise. $\mathcal{D}:=\diag{0,\ldots,n}$. $\mathcal{B}:=[0,\ldots,0,1]^T\in\mathbb{R}^{(n+1)}$.

\section{Problem statement and preliminaries}
\label{Sec:Prelim}

\subsection{Problem statement}

Consider as admissible signals $\mathcal{Y}_{(L,M)}^{(n+1)}$  the set of all signals $y : \mathbb{R}_+ \to \mathbb{R}$ that can be differentiated $n$ times for which $
\left|y^{(n+1)}(t)\right|\leq L(t) \text{ for almost all } t\geq 0$. Here,  $L:\mathbb{R}_+\to\mathbb{R}_+$ is a function such that 
\begin{equation}
\label{eq:M}
\frac{1}{L(t)}\left|\frac{d{L}(t)}{dt}\right|\leq M, \quad \forall t \geq 0,
\end{equation}
for some known constant $M \geq 0$. Note that $M$ is a bound for the logarithmic derivative of $L(t)$.

We consider the following problem:
\begin{problem}
\label{Problem}
Let $y\in\mathcal{Y}_{(L,M)}^{(n+1)}$ and consider a user-defined time $T_c>0$. Given measurements of $y(t)$ and knowledge of $L(t)$, the problem consists in estimating  the functions $y^{(i)}(t)$, $i=0,\cdots,n$, for all time $t\geq T_c$. 
\end{problem}

Let $z_i(t)$ denote the estimate for $y^{(i)}(t)$. Define  $z(t)=[z_0(t),\dots,z_n(t)]^T$.  In this paper, we will consider differentiators of the form 
\begin{align}
\dot{z}=-\mathcal{H}(z_0 - y ,t;T_c)+\mathcal{U}z, \label{Eq:Diff} 
\end{align}
where $\mathcal{U}$ is the matrix defined in the Notation and $\mathcal{H}:\mathbb{R}\times \mathbb{R}_+\to \mathbb{R}^{n+1}$ are correction functions to be designed. The correction functions have $T_c$ as a parameter. We restrict these functions to be continuous in $z_0 - y$ except at $z_0 = y$, and continuous in $t$ almost everywhere. Therefore, due to discontinuities, solutions to \eqref{Eq:Diff} are understood in  Filippov's sense \cite[Page 85]{Filippov1988DifferentialSides}.
 
To study the convergence of the differentiatior, consider the  differentiation error $e_i(t) = z_i(t) - y^{(i)}(t)$ for $i = 0, \cdots, n$. Its dynamics is given by 
\begin{align}
\dot{e}(t)=-\mathcal{H}(e_0(t),t;T_c)+\mathcal{U}e(t)-\mathcal{B}y^{(n+1)}(t), \label{Eq:DiffErr}
\end{align}
where $e=[e_0,\cdots,e_n]^T$ and $\mathcal{B}$ is the matrix defined in the Notation. Note that $y^{(n+1)}$ acts like a perturbation to the error dynamics, and that   $|y^{(n+1)}(t)|\leq L(t)$ by assumption.

To introduce the notion of settling-time, assume that $\mathcal{H}(e_0,t;T_c)$ is such that the origin of~\eqref{Eq:DiffErr} is asymptotically stable and  that it has unique solutions in forward-time for all $t\in[0,\infty)$. Then, the \textit{settling-time function} $T(e(0))$ of system~\eqref{Eq:DiffErr} for initial state $e(0) \in \mathbb R^{n+1}$ is
\begin{multline*}
T(e(0))=\\\inf\bigg\{\xi\geq 0: 
\forall y\in\mathcal{Y}_{(L,M)}^{(n+1)} , e(\xi';e(0),y)=0, \forall \xi'\geq \xi\bigg\},
\end{multline*}
where $e(t;e(0),y)$ is the solution of Eq. \eqref{Eq:DiffErr} for $t\geq 0$ with signal $y(t)$ and initial condition $e(0)$. We say that the origin of system~\eqref{Eq:DiffErr} is \textit{finite-time stable} if it is asymptotically stable~\citep{Khalil2002NonlinearSystems} and for every initial state $e(0) \in \mathbb R^n$, the settling-time function $T(e(0))$ is finite. We say that the origin of system~\eqref{Eq:DiffErr} is  \textit{fixed-time stable} if it is asymptotically stable~\citep{Khalil2002NonlinearSystems} and there exists $T_{{\max}}<\infty$ such that $T(e(0))\leq T_{\max}$ for all $e(0)\in\mathbb{R}^n$. Here,  $T_{\max}$ is the UBST of the system~\eqref{Eq:DiffErr}.

With the above definitions, a differentiator is \emph{exact} if the origin of its error dynamic is globally finite-time stable. Algorithm~\eqref{Eq:Diff} is said to be a differentiator with a predefined UBST if the origin of its differentiation error dynamic~\eqref{Eq:DiffErr} is fixed-time stable with a predefined UBST.

\subsection{Preliminaries}
\label{Appendix:Levant}

Here, we recall the design of Levant's finite-time differentiator with TVGs, which will be the base for our algorithm.

\begin{theorem}\citep{Levant2018GloballyGains}
\label{Th:Levant}
Consider the differentiator
\begin{align}
\dot{z}=&-\Phi(z_0 - y ,t;M,L(t))+\mathcal{U}z, \label{Eq:LevDiff} 
\end{align}
where 
\begin{multline*}
\Phi(w,t;M,L(t)):=[\phi_0(w,t;M,L(t)),\\ 
\dots,\phi_n(w,t;M,L(t))]^T
\end{multline*}
are functions recursively defined as  
$$ \phi_{i}(w,t;M,L(t)):=\chi_i( \ \phi_{i-1}(w,t;M,L(t)); \ M,L(t) ),$$
$i = 0, \cdots, n$, where $\phi_0(w,t;M,L(t)):=\chi_0(w;M,L(t))$ and 
\begin{multline}
\label{Eq:Levant1}
\chi_i(w;M,L(t)):=\lambda_{n-i}L(t)^{\frac{1}{n-i+1}}\sgn{w}^{\frac{n-i}{n-i+1}}+\mu_{n-i}Mw,
\end{multline}
with $\lambda_i$ and $\mu_i$ constant parameters. Then, there exist positive constants $\lambda_i$ and  $\mu_i$, $i=0, \ldots, n$, such that the algorithm of Eq.~\eqref{Eq:LevDiff} is an exact differentiator for any signal in $y \in \mathcal{Y}_{(L,M)}^{(n+1)}$. 
\end{theorem}
The above theorem implies that  the origin of the system 
\begin{align}
\dot{e}(t)=&-\Phi(e_0 (t),t;M,L(t))+\mathcal{U}e(t)-\mathcal{B}y^{(n+1)}(t), \label{Eq:LevDiffErr} 
\end{align}
is finite-time stable. In other words, the system~\eqref{Eq:DiffErr} is finite-time stable if we choose the correction functions as $\mathcal{H}(w,t;T_c)=\Phi(w,t;M,L(t))$.
Note  that finite-time stability is maintained if we relax the condition of Eq.~\eqref{eq:M} to: (i) $M$ is such that there exists $T^*\in[0,\infty)$ such that $L^{-1}(t)\left|d{L}(t)/dt\right|\leq M \text{ for all } t \geq T^*$; and  (ii) $L(t)$ is such that the solution of \eqref{Eq:DiffErr} exists for all $t\geq 0$.

\section{Main result}\label{Sec:ProbStatement}
\label{Sec:Main}

To solve  Problem~\ref{Problem} we ``redesign'' the correction functions   $\Phi(e_0,t;\bullet,\bullet)$ by combining them with a TVG $\kappa(\bullet)$ to obtain a new correction function $\mathcal{H}(e_0,t;T_c)$ such that Eq.~\eqref{Eq:Diff} is an exact differentiator and its UBST is predefined by $T_c$. We characterize $\kappa(\bullet)$ via an auxiliary function $\Omega(\bullet)$ as
$$\kappa(t) = \frac{d \varphi(t) }{d t}, $$
where $\varphi$ is defined using its inverse function  $\varphi^{-1}(\tau) := T_c \int_0^\tau \Omega(\xi) d \xi$. In this form, the class of TVG that our methodology considers is characterized by all functions $\Omega(\bullet)$ satisfying the following assumption:

\begin{assumption}
\label{Assump:NonAut}
The function $\Omega:\mathbb{R}_+\to\Bar{\mathbb{R}}_+\setminus\{0\}$ is such that: $\int_0^{\infty} \Omega(z)dz = 1$; $\Omega(\tau)<\infty$, for all $\tau>0$; it is either non-increasing or locally Lipschitz on $\mathbb{R}_+\setminus\{0\}$; and $\Omega(z)^{-1}\frac{d\Omega(z)}{dz}$ is uniformly bounded with respect to time and
satisfies
\begin{equation}
    \label{Eq:LimRho}
    \lim_{z\to\infty}\Omega(z)^{-1}\frac{d\Omega(z)}{dz}=-c,
\end{equation}
for some finite constant $c \geq 0$. 
\end{assumption}

Note that Assumption \ref{Assump:NonAut} implies $\varphi^{-1}(\tau)\leq T_c, \forall \tau \in[0,\infty)$. %
Table~\ref{tab:TVG} shows examples of TVG satisfying this assumption and the corresponding value of the $c$ constant.

Given a desired convergence time $T_c$, consider the redesigned correction function
\begin{equation}
 \mathcal{H}(e_0,t;T_c)=
\left\lbrace
\begin{array}{cl}
  \Gamma(e_0,t;M,L(t) )   & \text{ for } t\in[0,T_c),\\
  \Phi(e_0,t;M,L(t))   & \text{ otherwise,}
\end{array}
\right.
\label{eq:H}
\end{equation}
with $\Phi(e_0, t; M, L(t))$ chosen as in Theorem~\ref{Th:Levant} and 
\begin{equation}
    \begin{split}
    \Gamma(e_0,t;M, L(t)) = \ \Lambda(t) \big[ & \mathcal{Q}(c)\Phi(e_0;\mathcal{M},L(t)\kappa(t)^{-(n+1)})   \\
    & + (\mathcal{U}-c \mathcal{D})^{n+1} \mathcal{B}e_0 \big].
    \end{split}
\end{equation} 
Here, $\Lambda(t)=\diag{\kappa(t),\cdots,\kappa(t)^{n+1}}$ is a gain matrix built from the scalar TVG, $\mathcal{M}>(n+1)c$ with $c$ the TVG constant of Assumption 1, and $\mathcal{Q}(c): = \begin{bmatrix}
(\mathcal{U}-c \mathcal{D})^{n} \mathcal{B}; &
\cdots; &
(\mathcal{U}-c \mathcal{D}) \mathcal{B}; &
\mathcal{B}
\end{bmatrix}$
with the matrices $\mathcal{B}$, $\mathcal{D}$  and $\mathcal{U}$  given in the Notation. The intuition about Eq. \eqref{eq:H} is as follows. Before $T_c$, the correction function consists of the original function of Levant's differentiator redesigned with the TVG $\Lambda(t)$, a rescaling $\mathcal Q(c)$, and a linear term. Such redesign ensures its convergence before  $T_c$. In particular, $\Lambda(t)$ acts as a time scaling providing such predefined-time convergence.  After $T_c$, the correction function is switched to Levant's differentiator, maintaining exactness. More precisely, our main result is the following:

\begin{theorem}
\label{Theorem:Main} Consider the differentiator algorithm ~\eqref{Eq:Diff} with the correction function $\mathcal H(e_0, t; T_c)$ of Eq. \eqref{eq:H}. Then, for any $e(0)\in\mathbb{R}^{n+1}$, the algorithm has a unique Filippov solution\footnote{Note that the standard Filippov conditions for the existence of solutions \cite[Page 85, Theorem 8]{Filippov1988DifferentialSides} are not satisfied at $t=T_c$ for algorithm~\eqref{Eq:Diff}. However, the existence and uniqueness of Filippov solutions are guaranteed for the proposed design through a time-scaling argument.} defined for all $t\geq 0$, and it solves Problem~\ref{Problem} (i.e., the origin of the differentiation error dynamics is fixed-time stable with $T_c$ as a predefined UBST). Moreover, if $L(t)$ is such that for the base differentiator of Eq.~\eqref{Eq:LevDiff}, the settling-time function $\mathcal{T}$ of its differentiation error dynamics satisfies 
\begin{equation}
\label{Eq:Sumpreme}
\sup_{e(0)\in\mathbb{R}^{n+1}}\mathcal{T}(e(0))=\infty,
\end{equation} 
then $T_c$ is the least UBST.%
\end{theorem}

\begin{table}
    \centering
    \begin{tabular}{|c|c|c|}
    \hline
      & time-varying gain  $\kappa(t)$  & constant  $c$\\
      \hline\hline
      \textit{(i)} & $\alpha^{-1}(T_c-t)^{-1}$ &  $\alpha$\\
      \hline
      \textit{(ii)} &  $\frac{\pi}{2}\sec(\frac{\pi t}{2T_c})^2$ & $0$\\
      \hline
      \textit{(iii)} & $\frac{\gamma}{T_c} \tan(\gamma \frac{t}{T_c}+\frac{\pi}{2}-\gamma)$  & $1$\\
       
      \hline
      \textit{(iv)} & $\frac{t+\beta}{\alpha(T_c-t)}$ & $\frac{\alpha}{1+\beta}$\\
      \hline
    \end{tabular}
    \caption{Examples of TVGs satisfying Assumption 1. Here, $\alpha > 0,\beta>0$, and $0 < \gamma < \pi/2$.}
    \label{tab:TVG}
\end{table}

\begin{remark} 
\label{Rem:Filtering}
Following~\cite{Levant2019,Carvajal-Rubio2022}, we can extend our methodology to redesign filtering differentiators by using $n=n_d+n_f$ and considering signals in $\mathcal{Y}^{(n_d+1)}_{(L,M)}$, where $n_f$ is the filtering order, as their filtering properties can be very useful in the presence of measurement noise~\citep{Levant2019}.
\end{remark}

\section{Obtaining a bounded gain and accuracy to measurement noise}
\label{sec:practical}

In the proof of Theorem~\ref{Theorem:Main}, we show that the settling time of the differentiator error actually satisfies $T(e(0))<T_c$, implying that the differentiator always converges \emph{before} $T_c$. Therefore, since the time-varying gain grows unbounded only at $T_c$,  $\kappa(t)$ will be finite at the settling time $T(e(0))$. However, as the initial differentiation error $e(0)$ grows, the settling time $T(e(0))$ will approach $T_c$, i.e.,
\begin{equation}
\label{Eq:SumpremeST}
\sup_{e(0)\in\mathbb{R}^{n+1}}T(e(0))=T_c,
\end{equation}
and the TVG will grow unbounded at $T_c$, making the gain $\kappa(T(e(0)))$ at the settling time  grow unbounded as a function of the initial differentiation error. Such unbounded gain is a drawback in our methodology, making its application challenging in the presence of measurement noise. In particular, arbitrarily small noise may cause divergence of solutions at $t=T_c$. Furthermore, an arbitrary small Lipschitz measurement noise can make the differentiation error arbitrarily large at $T_c$ ~\citep{Aldana-Lopez2022OnAlgorithms}.

Our methodology offers a simple workaround to circumvent the above challenges: instead of switching from $\Gamma(e_0, t; M, L(t))$ to $\Phi(e_0, t; M, L(t))$ at time $T_c$ in the correction function of Eq. \eqref{eq:H}, switch at an earlier time $T_c^*  < T_c$. With this workaround, the differentiation error of our algorithm remains finite-time convergent. Furthermore, there exists a compact neighborhood $R \subset \mathbb R^{n+1}$ of initial differentiation errors around the origin whose settling time is bounded by $T_c$. Such neighborhood size can be arbitrarily large by choosing $\alpha$ and $\mathcal{M}$ adequately. In addition, when the signal to be differentiated is measured with noise $\eta_0(t)$,  the accuracy of the redesigned differentiator remains similar to the original differentiator in the following sense:

\begin{proposition}
\label{prop:noise}
Consider the differentiator algorithm  \eqref{Eq:Diff} with the correction function of Eq. \eqref{eq:H} applied to the signal $y = y_0 + \eta_0$, where $y_0 \in \mathcal Y^{(n+1)}_{L, M}$ is the nominal part of the signal and $\eta_0: [0, \infty) \rightarrow \mathbb R]$ is a measurable function representing measurement noise. Suppose that $|\eta_0(t)|\leq \eta L(t)$ for some constant $\eta>0$ and that the switching time  in Eq. \eqref{eq:H} is chosen as $T_c^* < T_c$.  Then, there exist  $\eta$ and $s$ small enough and a $T_c^*$ sufficiently close to $T_c$ such that
\begin{equation}
    |e_i(t)|\leq  \tilde{\gamma}_iL(t)\frac{\kappa(T_c^*)^{n+1}}{\kappa(0)^{n-i+1}}\eta^{\frac{n-i+1}{n+1}}, i=0,\dots, n,
\end{equation}
for some $\tilde{\gamma}_0,\dots,\tilde{\gamma}_n>0$ and $t\in[s,T_c^*]$.
\end{proposition}

\section{Numerical examples and comparisons}
\label{Sec:Discussion}

Here we provide numerical examples to illustrate the advantages and limitations of our results. All simulations were created in OpenModelica using the Euler integration method with a time step $t_{\textrm{step}}=2\times10^{-4}$. These simulations use the workaround suggested in Section~\ref{sec:practical} to maintain the TVG bounded thus avoiding numerical issues during the simulation even without measurement noise. In simulations, applying the workaround suggested in Section~\ref{sec:practical} requires choosing $T_c^* \le T_c - t_\textrm{step}$. If $T_c$ is a multiple of the integration time step, then $T_c - t_\textrm{step}$ is the last integration step before the singularity of the TVG. Importantly, recall that Proposition~\ref{prop:noise} indicates that the sensitivity to noise increases as $T_c^*$ approaches $T_c$. Following the original results of \citep{Levant2018GloballyGains}, to build the correction function $\Phi(e_0, t; M, L(t))$ of Theorem 1  we choose its parameters  as
\begin{multline}
\{\lambda_i,\mu_i\}_{i=0}^{n}=(1.1,2),(1.5,3),(2,4),(3,7),(5,9),\ldots
\end{multline}

\subsection{High-order exact differentiators with a predefined UBST}

We start by discussing a second-order differentiator for signals with exponential growth.

\begin{example}
\label{Ex:SOD} 
Consider the differentiator algorithm~\eqref{Eq:Diff} with $n=2$ and let  $y_0(t)=2\sin\left(\frac{1}{2} t^2\right)$. Note that this signal satisfies
$$|{y}^{(3)}(t)|\leq L(t),$$ with $L(t)=\sqrt{1+4 t^6+36 t^2}$. Here,   $\frac{1}{L(t)}\left|\frac{d {L}(t)}{dt}\right|\leq M=3.5$. Consider our algorithm with $\kappa(t)$ given in Table~\ref{tab:TVG}-(iii) with $\gamma=0.01$. Note that this TVG has a constant $c =1$. Accordingly, we choose $\mathcal{M}=8 > c(n+1)$ and
\begin{equation}
\label{Eq:QExa1}
\mathcal{Q}(c)=
\begin{bmatrix}
     1  &   0  &   0 \\
    -3  &   1  &   0 \\
     4  &  -2  &   1 
\end{bmatrix} \text{ and } \ 
(\mathcal{U}-c\mathcal{D})^3\mathcal{B}=
\begin{bmatrix}
-3 \\ 7 \\-8
\end{bmatrix}.
\end{equation}
Hence, the proposed differentiator for $t\in[0,T_c^*)$ takes the form
\begin{align}
\dot{z}_0=&-\kappa(t)(\varrho_0(e_0,t)-3e_0)+z_1,\\
\dot{z}_1=&-\kappa(t)^2(-3\varrho_0(e_0,t)+\varrho_1(e_0,t)+7e_0)+z_2, \label{Eq:SODiff}\\
\dot{z}_2=&-\kappa(t)^3(4\varrho_0(e_0,t)-2\varrho_1(e_0,t)+\varrho_2(e_0,t)-8e_0),
\end{align}
where
$
\varrho_i(e_0,t):=\phi_i(e_0,t;\mathcal{M},L(t)\kappa(t)^{-3})
$
with $\{\phi_i\}_{i=0}^2$ the error correction functions of the time-varying differentiator of Theorem~\ref{Th:Levant}.

The convergence of the above algorithm for different initial conditions is shown in the first row of Fig.~\ref{Fig:SOD}, where $T_c=5$. To maintain the gain $\kappa(t)$ bounded, we select $T_c^*=4.5$. It can be observed that exact convergence before the desired time is obtained even for large initial conditions.

\begin{figure*}
    \centering
    \def\svgwidth{16cm}
\begingroup%
  \makeatletter%
  \providecommand\color[2][]{%
    \errmessage{(Inkscape) Color is used for the text in Inkscape, but the package 'color.sty' is not loaded}%
    \renewcommand\color[2][]{}%
  }%
  \providecommand\transparent[1]{%
    \errmessage{(Inkscape) Transparency is used (non-zero) for the text in Inkscape, but the package 'transparent.sty' is not loaded}%
    \renewcommand\transparent[1]{}%
  }%
  \providecommand\rotatebox[2]{#2}%
  \newcommand*\fsize{\dimexpr\f@size pt\relax}%
  \newcommand*\lineheight[1]{\fontsize{\fsize}{#1\fsize}\selectfont}%
  \ifx\svgwidth\undefined%
    \setlength{\unitlength}{1080bp}%
    \ifx\svgscale\undefined%
      \relax%
    \else%
      \setlength{\unitlength}{\unitlength * \real{\svgscale}}%
    \fi%
  \else%
    \setlength{\unitlength}{\svgwidth}%
  \fi%
  \global\let\svgwidth\undefined%
  \global\let\svgscale\undefined%
  \makeatother%
  \begin{picture}(1,0.43)%
    \lineheight{1}%
    \setlength\tabcolsep{0pt}%
    \put(0,0){\includegraphics[width=\unitlength,page=1]{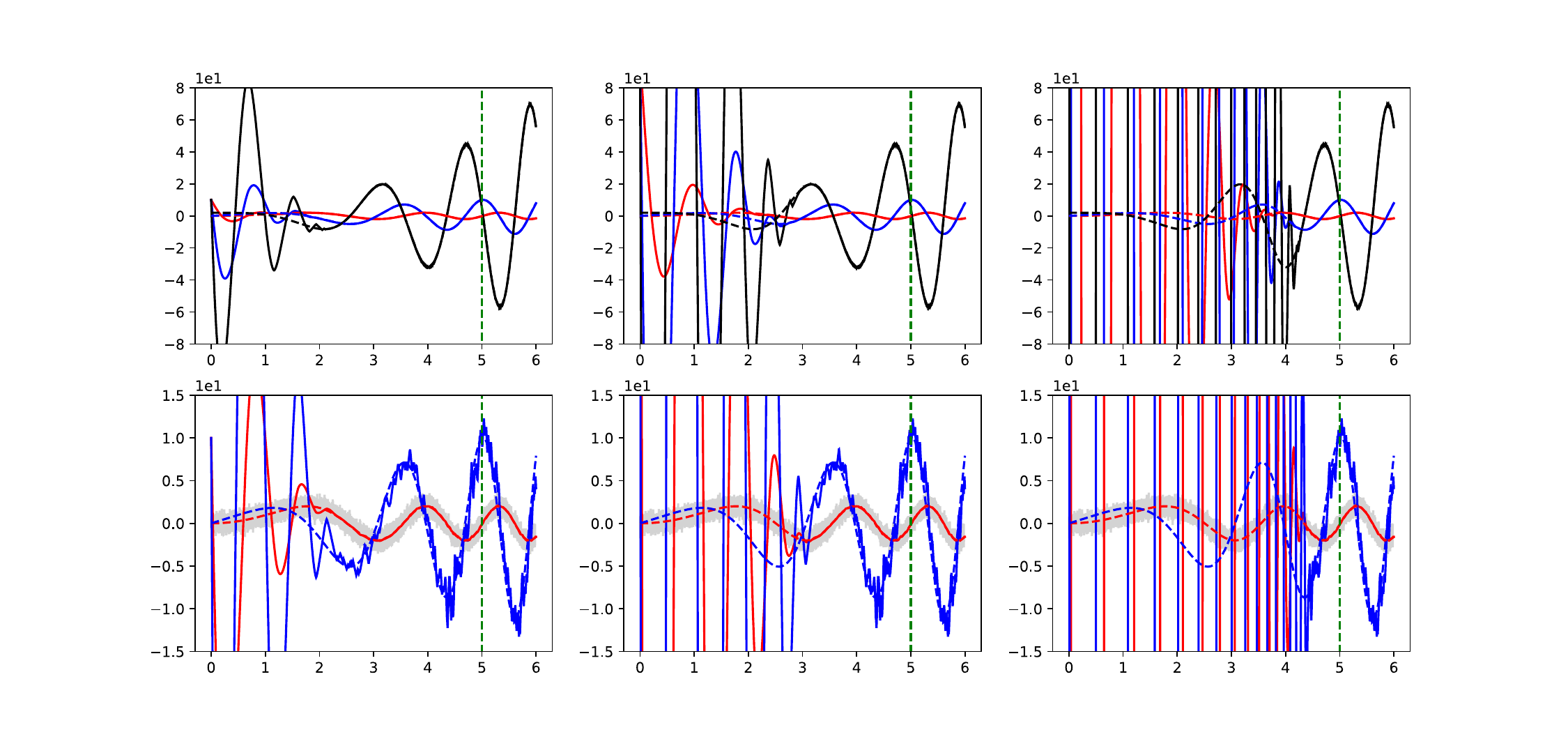}}%
    \put(0,0){\includegraphics[width=\unitlength,page=1]{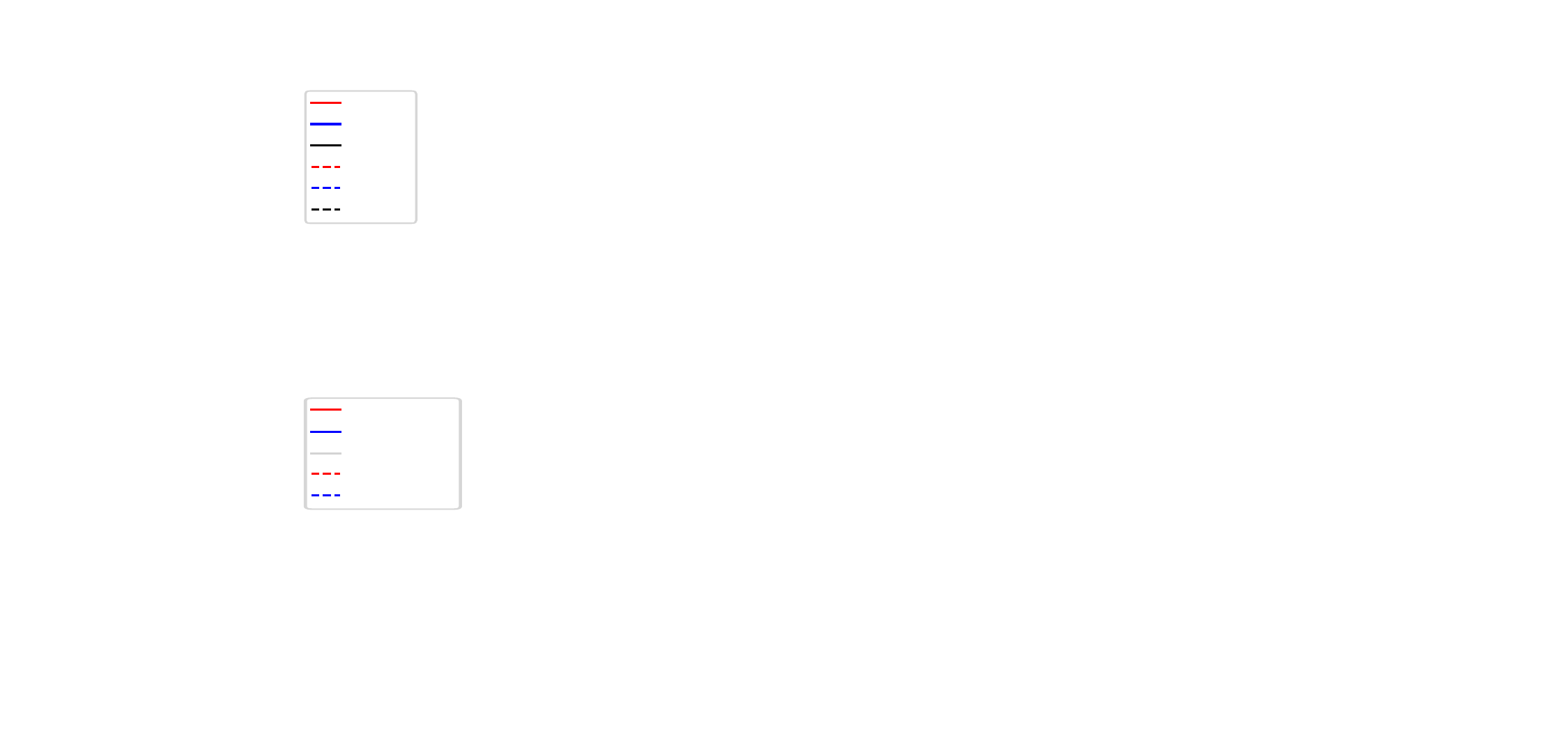}}%
    \tiny{
    \put(0.22067925,0.02442618){\color[rgb]{0,0,0}\makebox(0,0)[lt]{\lineheight{1.25}\smash{\begin{tabular}[t]{l}time\end{tabular}}}}%
    \put(0.49845705,0.02442618){\color[rgb]{0,0,0}\makebox(0,0)[lt]{\lineheight{1.25}\smash{\begin{tabular}[t]{l}time\end{tabular}}}}%
    \put(0.77623483,0.02442618){\color[rgb]{0,0,0}\makebox(0,0)[lt]{\lineheight{1.25}\smash{\begin{tabular}[t]{l}time\end{tabular}}}}%
    \put(0.09002649,0.0544827){\color[rgb]{0,0,0}\rotatebox{90}{\makebox(0,0)[lt]{\lineheight{1.25}\smash{\begin{tabular}[t]{l}$n_f=1$, with Noise\end{tabular}}}}}%
    \put(0.09002649,0.25031603){\color[rgb]{0,0,0}\rotatebox{90}{\makebox(0,0)[lt]{\lineheight{1.25}\smash{\begin{tabular}[t]{l}$n_f=0 $, without Noise\end{tabular}}}}}%
    \put(0.15796366,0.41716886){\color[rgb]{0,0,0}\makebox(0,0)[lt]{\lineheight{1.25}\smash{\begin{tabular}[t]{l}$z_0(t)=z_1(0)=z_2(0)=1e1$\end{tabular}}}}%
    \put(0.43574145,0.41716886){\color[rgb]{0,0,0}\makebox(0,0)[lt]{\lineheight{1.25}\smash{\begin{tabular}[t]{l}$z_0(t)=z_1(0)=z_2(0)=1e2$\end{tabular}}}}%
    \put(0.71351923,0.41716886){\color[rgb]{0,0,0}\makebox(0,0)[lt]{\lineheight{1.25}\smash{\begin{tabular}[t]{l}$z_0(t)=z_1(0)=z_2(0)=1e4$\end{tabular}}}}%
    \put(0.15796366,0.21994663){\color[rgb]{0,0,0}\makebox(0,0)[lt]{\lineheight{1.25}\smash{\begin{tabular}[t]{l}$w_1(0)=z_0(0)=z_1(0)=1e1$\end{tabular}}}}%
    \put(0.43574145,0.21994663){\color[rgb]{0,0,0}\makebox(0,0)[lt]{\lineheight{1.25}\smash{\begin{tabular}[t]{l}$w_1(0)=z_0(0)=z_1(0)=1e2$\end{tabular}}}}%
    \put(0.71351923,0.21994663){\color[rgb]{0,0,0}\makebox(0,0)[lt]{\lineheight{1.25}\smash{\begin{tabular}[t]{l}$w_1(0)=z_0(0)=z_1(0)=1e4$\end{tabular}}}}%
    \put(0.22362802,0.39744799){\color[rgb]{0,0,0}\makebox(0,0)[lt]{\lineheight{1.25}\smash{\begin{tabular}[t]{l}$z_0(t)$\end{tabular}}}}%
    \put(0.22362802,0.38494799){\color[rgb]{0,0,0}\makebox(0,0)[lt]{\lineheight{1.25}\smash{\begin{tabular}[t]{l}$z_1(t)$\end{tabular}}}}%
    \put(0.22362802,0.37244799){\color[rgb]{0,0,0}\makebox(0,0)[lt]{\lineheight{1.25}\smash{\begin{tabular}[t]{l}$z_2(t)$\end{tabular}}}}%
    \put(0.22362802,0.35994799){\color[rgb]{0,0,0}\makebox(0,0)[lt]{\lineheight{1.25}\smash{\begin{tabular}[t]{l}$y(t)$\end{tabular}}}}%
    \put(0.22362802,0.34605911){\color[rgb]{0,0,0}\makebox(0,0)[lt]{\lineheight{1.25}\smash{\begin{tabular}[t]{l}$\dot{y}(t)$\end{tabular}}}}%
    \put(0.22362802,0.33217022){\color[rgb]{0,0,0}\makebox(0,0)[lt]{\lineheight{1.25}\smash{\begin{tabular}[t]{l}$\ddot{y}(t)$\end{tabular}}}}%
    \put(0.22362802,0.18911466){\color[rgb]{0,0,0}\makebox(0,0)[lt]{\lineheight{1.25}\smash{\begin{tabular}[t]{l}$z_1(t)$\end{tabular}}}}%
    \put(0.22362802,0.17661466){\color[rgb]{0,0,0}\makebox(0,0)[lt]{\lineheight{1.25}\smash{\begin{tabular}[t]{l}$y(t)+\eta_0(t)$\end{tabular}}}}%
    \put(0.22362802,0.16272577){\color[rgb]{0,0,0}\makebox(0,0)[lt]{\lineheight{1.25}\smash{\begin{tabular}[t]{l}$y(t)$\end{tabular}}}}%
    \put(0.22362802,0.14883688){\color[rgb]{0,0,0}\makebox(0,0)[lt]{\lineheight{1.25}\smash{\begin{tabular}[t]{l}$\dot{y}(t)$\end{tabular}}}}%
    \put(0.22362802,0.20439244){\color[rgb]{0,0,0}\makebox(0,0)[lt]{\lineheight{1.25}\smash{\begin{tabular}[t]{l}$z_0(t)$\end{tabular}}}}%
    }
  \end{picture}%
\endgroup%
    \caption{In the first row, simulation of Example~\ref{Ex:SOD}, with desired UBST given by $T_c=5$, using our algorithm~\eqref{Eq:Diff} with $n=2$. In the second row, simulation of Example~\ref{Ex:FilteringDiff}, using our algorithm in the form of a filtering differentiator~\eqref{Eq:ExFiltering} with $n_f=1$, $n_d=1$ and the signal $y$ contaminated with a white Gaussian noise signal $\eta_0(t)$ with zero mean and standard deviation $\sigma=0.5$.}
    \label{Fig:SOD}
\end{figure*}
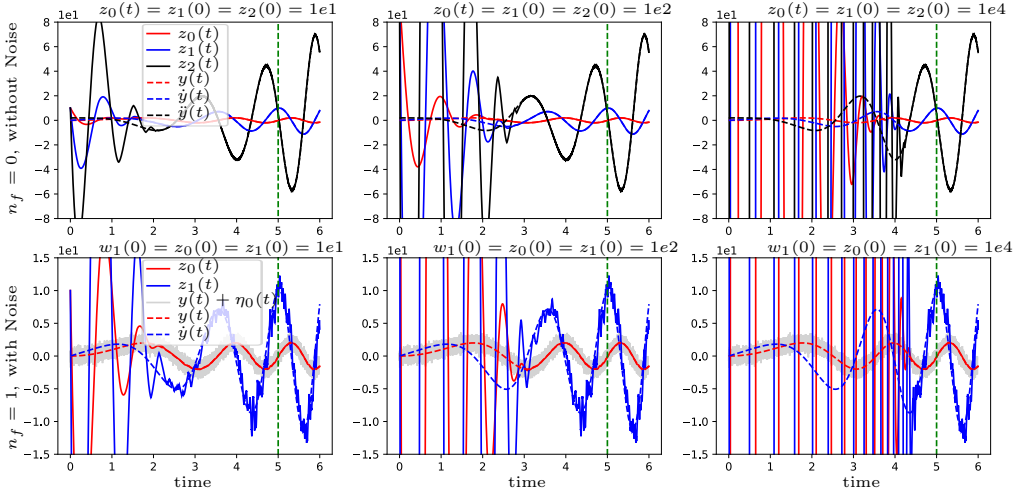
\end{example}

The next example  shows our algorithm in the form of a filtering differentiator as discussed in Remark~\ref{Rem:Filtering} with $n=2$, $n_d=1$ and $n_f=1$. Here, we consider  signals having a second-order derivative bounded by $L(t)=\sqrt{4+4t^4}$.

\begin{example}
\label{Ex:FilteringDiff}
Consider the signal to be differentiated to be contaminated by Gaussian measurement noise as $y(t)=y_0(t)+\eta_0(t)$ where $\eta_0(t)$ has zero mean and standard deviation $\sigma=0.5$ with the nominal signal  $y_0(t)=2\sin\left(\frac{1}{2} t^2\right)$ satisfying $|\ddot{y}_0(t)|\leq L(t)$ with $L(t)=\sqrt{4+4t^4}$. Note that  $\frac{1}{L(t)}\left|\frac{d {L}(t)}{dt}\right|\leq M=3.5$. Consider our algorithm with $\kappa(t)$ given in Table~\ref{tab:TVG}-(iii) with $\gamma=0.01$ and let $\mathcal{M}=8$. To maintain  $\kappa(t)$ bounded, we select $T_c^*=4.5$. Thus, $n=2$ and $c=1$ as in Example~\ref{Ex:SOD}, then $\mathcal{Q}(c)$ and $(\mathcal{U}-c\mathcal{D})^3\mathcal{B}$ are given as in~\eqref{Eq:QExa1}, and the filtering differentiator for $t\in[0,T_c^*)$ is given by 
\begin{align}
\dot{w}_1=&-\kappa(t)(\varrho_0(w_1,t)-3w_1)+z_0-y,\\
\dot{z}_0=&-\kappa(t)^2(-3\varrho_0(w_1,t)+\varrho_1(w_1,t)+7w_1)+z_1, \label{Eq:ExFiltering}\\
\dot{z}_1=&-\kappa(t)^3(4\varrho_0(w_1,t)-2\varrho_1(w_1,t)+\varrho_2(w_1,t)-8w_1),
\end{align}
where
$
\varrho_i(e_0,t):=\phi_i(e_0,t;\mathcal{M},L(t)\kappa(t)^{-3})
$
with $\{\phi_i\}_{i=0}^2$ the error correction functions of the time-varying differentiator of Theorem~\ref{Th:Levant}. The convergence of this algorithm for different initial conditions is shown in the second row of Fig.~\ref{Fig:SOD}.
\end{example}

\subsection{Comparison with arbitrary order exact differentiators for polynomial signals of $n-$th order}

Higher-order exact differentiators with a predefined UBST are  provided in~\cite{Menard2017,Holloway2019}, but they are designed for polynomial signals of $n-$th order only. Whereas the UBST in~\cite{Menard2017} is very conservative, the convergence in the algorithm given by~\cite{Holloway2019} occurs precisely at the predefined-time. 

The TVG of Table~\ref{tab:TVG}-\textit{(i)} was also used in~\cite{Holloway2019} with $\alpha=1$. However, in such an algorithm, $\lim_{t\to T(e(0))} \kappa(t)=\infty$ for every nonzero $e(0)$. The workaround discussed in Section~\ref{sec:practical} to maintain the TVG bounded was also suggested by~\cite{Holloway2019}. However, with such workaround an exact differentiator is no longer obtained with the algorithm by~\cite{Holloway2019} (recall that convergence occurs precisely at $T_c$). In fact, the magnitude of the error at $T_c^*$, $\|e(T_c^*)\|$, grows with the initial condition. We illustrate this case in the following example.

\begin{example}
\label{Ex:Krstic} 
Consider the second-order differentiator problem. To provide a fair comparison to~\cite{Holloway2019}, which we consider the closest to our approach, consider the quadratic polynomial function $y(t)=\frac{1}{2}t^2+t+1$.  Notice that the algorithm by~\cite{Holloway2019} is given by Eq.~\eqref{Eq:Diff} with $h_i(e_0,t;T_c)=g_i(e_0,t;T_c)e_0,$ and
\begin{align}
    g_0(e_0,t;T_c)=&l_1+T_c\kappa(t)\left(\frac{p_{11}(m+3)}{T_c}-p_{21}\right)\\
    g_1(e_0,t;T_c)=&l_2+T_c^2\kappa(t)^2\left(\frac{p_{21}(m + 4)}{T_c}-p_{31}\right)\\& -  p_{21}T_c\kappa(t)g_0(e_0,t;T_c) \label{Eq:ObsHolloway}\\
    g_2(e_0,t;T_c)=&l_3+T_c^3\kappa(t)^3\frac{p_{31}(m+5)}{T_c} \\&-  p_{31}T_c^2\kappa(t)^2g_0(e_0,t;T_c)  \\&-p_{32}T_c\kappa(t)g_1(e_0,t;T_c)
\end{align}
where $p_{11}=1, p_{21}=-\frac{2(m+3)}{T_c}, p_{31}=\frac{(m+3)(m+4)}{T_c^2}, p_{32}=-\frac{m+3}{T_c}$, $m$ is a design parameter, here chosen to be $m=1$, and $[l_1,l_2,l_3]=[6,11,6]$. Moreover, $\kappa(t)$ is given in Table~\ref{tab:TVG}-(i) with $\alpha=1$. Consider the workaround discussed in Section~\ref{sec:practical} to maintain a bounded TVG. For illustration purposes, we choose $\kappa(t)\leq 10$. Thus, $T_c^*=0.9$. Notice that, since the error system with the algorithm~\eqref{Eq:ObsHolloway} is linear time-varying with bounded systems matrix for $t\in[0,T_c^*]$, then zero differentiation error is no longer obtained. Furthermore, the magnitude $\|e(T_c^*;e(0),y)\|$ grows linearly with the initial condition $e(0)$.%

For our algorithm, consider $\kappa(t)$ given by Table~\ref{tab:TVG}-(i) with $c=\alpha=1$, $L(t)=0.1\exp{-5t}$, $M=0.5$, and $\mathcal{M}=20$. Thus, $n=2$, as in Example~\ref{Ex:SOD}, then $\mathcal{Q}(c)$ and $(\mathcal{U}-c\mathcal{D})^3\mathcal{B}$ are given as in~\eqref{Eq:QExa1} and the differentiator structure is given as in~\eqref{Eq:SODiff}.
The convergence is shown in the second row of Fig.~\ref{Fig:Krstic}. Notice that, under the same workaround, even for an initial condition satisfying $z_0(0)=z_1(0)=z_2(0)= 10000$ zero differentiation error before $T_c$ is still obtained. A significant advantage in contrast with the algorithm by~\cite{Holloway2019}.

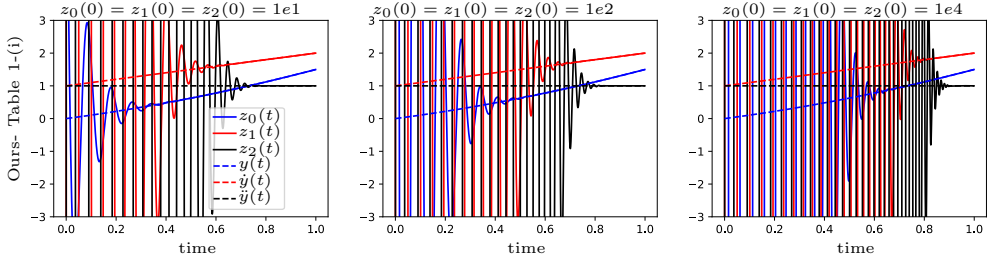
\begin{figure*}
    \centering
    \def\svgwidth{16cm}
\begingroup%
  \makeatletter%
  \providecommand\color[2][]{%
    \errmessage{(Inkscape) Color is used for the text in Inkscape, but the package 'color.sty' is not loaded}%
    \renewcommand\color[2][]{}%
  }%
  \providecommand\transparent[1]{%
    \errmessage{(Inkscape) Transparency is used (non-zero) for the text in Inkscape, but the package 'transparent.sty' is not loaded}%
    \renewcommand\transparent[1]{}%
  }%
  \providecommand\rotatebox[2]{#2}%
  \newcommand*\fsize{\dimexpr\f@size pt\relax}%
  \newcommand*\lineheight[1]{\fontsize{\fsize}{#1\fsize}\selectfont}%
  \ifx\svgwidth\undefined%
    \setlength{\unitlength}{1080bp}%
    \ifx\svgscale\undefined%
      \relax%
    \else%
      \setlength{\unitlength}{\unitlength * \real{\svgscale}}%
    \fi%
  \else%
    \setlength{\unitlength}{\svgwidth}%
  \fi%
  \global\let\svgwidth\undefined%
  \global\let\svgscale\undefined%
  \makeatother%
  \begin{picture}(1,0.24)%
    \lineheight{1}%
    \setlength\tabcolsep{0pt}%
    \put(0,0){\includegraphics[width=\unitlength,page=1,trim=0cm 0cm 0cm 9cm,clip]{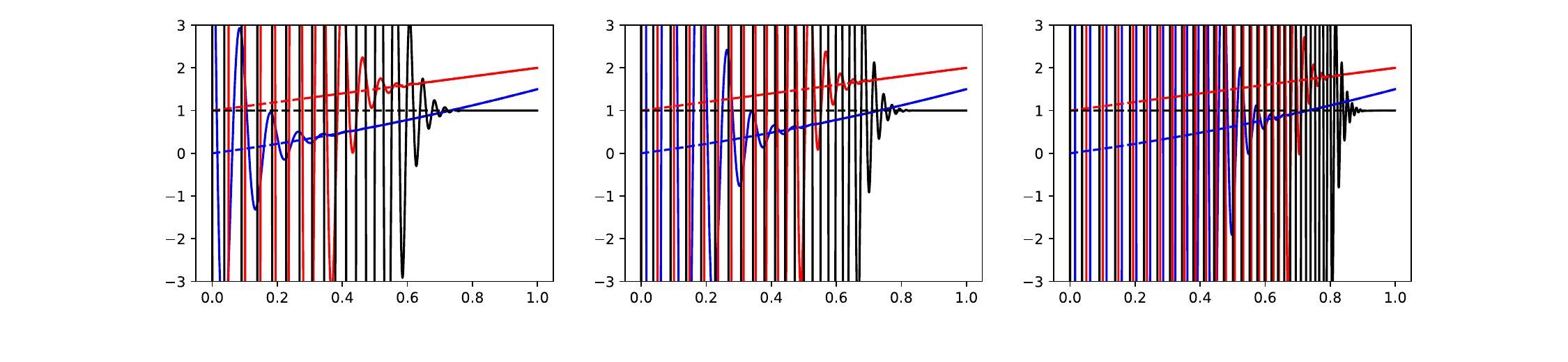}}%
    \put(0,0){\includegraphics[width=\unitlength,page=1,trim=0cm 0cm 0cm 9cm,clip]{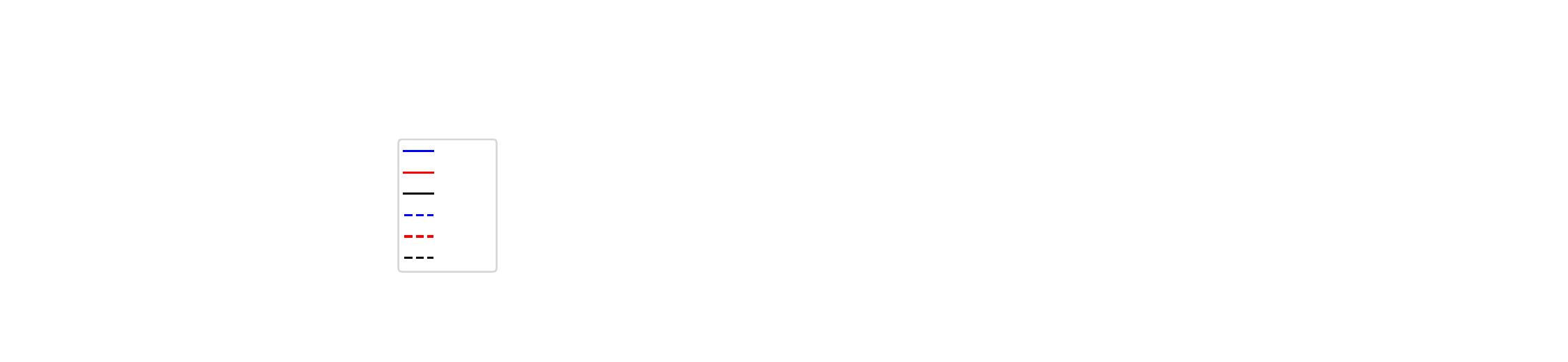}}%
    \tiny{
    \put(0.22863079,0.0212378){\color[rgb]{0,0,0}\makebox(0,0)[lt]{\lineheight{1.25}\smash{\begin{tabular}[t]{l}time\end{tabular}}}}%
    \put(0.77568958,0.0212378){\color[rgb]{0,0,0}\makebox(0,0)[lt]{\lineheight{1.25}\smash{\begin{tabular}[t]{l}time\end{tabular}}}}%
    \put(0.50216019,0.0212378){\color[rgb]{0,0,0}\makebox(0,0)[lt]{\lineheight{1.25}\smash{\begin{tabular}[t]{l}time\end{tabular}}}}%
    \put(0.27846996,0.13222108){\color[rgb]{0,0,0}\makebox(0,0)[lt]{\lineheight{1.25}\smash{\begin{tabular}[t]{l}$z_0(t)$\end{tabular}}}}%
    \put(0.27846996,0.11833219){\color[rgb]{0,0,0}\makebox(0,0)[lt]{\lineheight{1.25}\smash{\begin{tabular}[t]{l}$z_1(t)$\end{tabular}}}}%
    \put(0.27846996,0.1044433){\color[rgb]{0,0,0}\makebox(0,0)[lt]{\lineheight{1.25}\smash{\begin{tabular}[t]{l}$z_2(t)$\end{tabular}}}}%
    \put(0.27846996,0.0919433){\color[rgb]{0,0,0}\makebox(0,0)[lt]{\lineheight{1.25}\smash{\begin{tabular}[t]{l}$y(t)$\end{tabular}}}}%
    \put(0.27846996,0.07805441){\color[rgb]{0,0,0}\makebox(0,0)[lt]{\lineheight{1.25}\smash{\begin{tabular}[t]{l}$\dot{y}(t)$\end{tabular}}}}%
    \put(0.27846996,0.06416553){\color[rgb]{0,0,0}\makebox(0,0)[lt]{\lineheight{1.25}\smash{\begin{tabular}[t]{l}$\ddot{y}(t)$\end{tabular}}}}%
    \put(0.09625044,0.08){\color[rgb]{0,0,0}\rotatebox{90}{\makebox(0,0)[lt]{\lineheight{1.25}\smash{\begin{tabular}[t]{l}Ours- Table~\ref{tab:TVG}-(i)\end{tabular}}}}}%
    \put(0.13,0.22){\color[rgb]{0,0,0}\makebox(0,0)[lt]{\lineheight{1.25}\smash{\begin{tabular}[t]{l}$z_0(0)=z_1(0)=z_2(0)= 1e1$\end{tabular}}}}%
    \put(0.68,0.22){\color[rgb]{0,0,0}\makebox(0,0)[lt]{\lineheight{1.25}\smash{\begin{tabular}[t]{l}$z_0(0)=z_1(0)=z_2(0)= 1e4$\end{tabular}}}}%
    \put(0.39,0.22){\color[rgb]{0,0,0}\makebox(0,0)[lt]{\lineheight{1.25}\smash{\begin{tabular}[t]{l}$z_0(0)=z_1(0)=z_2(0)= 1e2$\end{tabular}}}}%
    }
  \end{picture}%
\endgroup%
    \caption{Simulation of Example~\ref{Ex:Krstic}, with desired UBST given by $T_c=1$.}
    \label{Fig:Krstic}
\end{figure*}

\end{example}

 \subsection{Comparison with autonomous first-order differentiator}

As mentioned in the introduction, predefined-time autonomous exact differentiators for signals whose $n-$th derivative is Lipschitz have only been proposed in the literature for the case where $n=1$ by~\cite{Cruz-Zavala2011,Seeber2021robust}. In the next example we compare it against the algorithm by~\cite{Seeber2021robust}, using the parameters proposed in \cite{Seeber2020ExactBound}.

\begin{example}
 \label{Ex:Seeber}

Let $y(t)=0.75cos(t)+0.0025sin(10t)+t$ with $L(t)=L=1$. For comparison, consider the algorithm in~\citep{Seeber2020ExactBound}, i.e., algorithm~\eqref{Eq:Diff} with $h_i(w,t;T_c)=k_i\nu_i(w;T_c)$, $i=1,2$, where
 \begin{align}
     \nu_1(w;T_c)&=\sgn{w}^{\frac{1}{2}}+k_3^2\sgn{w}^{\frac{3}{2}}\\
     \nu_2(w;T_c)&=\sgn{w}^{0}+4k_3^2w+3k_3^4\sgn{w}^{2}, \label{Eq:Seeber}
 \end{align}
 where $k_1=4\sqrt{L}$, $k_2=2L$ and $k_3=\frac{9.8}{T_c\sqrt{L}}$. 
 
 Now consider our algorithm with $\kappa(t)$ given in Table~\mbox{\ref{tab:TVG}-(i)} with $\alpha=1$, $L=1$, $M=0.1$ and $\mathcal{M}=3.2$. Thus, $c=1$ and the differentiator for $t\in[0,T
 ^*)$ becomes:
 \begin{align}
     \dot{z}_0&=-\kappa(t)(\varrho_0(e_0,t)-c e_0)+z_1\\
     \dot{z}_1&=-\kappa(t)(\varrho_1(e_0,t)-c\varrho_0(e_0,t)+c^2e_0)
 \end{align}
 where
 $
 \varrho_i(e_0,t):=\phi_i(e_0,t;\mathcal{M},L(t)\kappa(t)^{-2})
 $
 with $\{\phi_i\}_{i=0}^1$, the error correction functions of the time-varying differentiator of Theorem~\ref{Th:Levant}. For the first-order case, the error correction functions of Theorem~\ref{Th:Levant} can be written in a recursive form as~\cite{Levant2018GloballyGains}:
\begin{align}
\phi_1(w,t;M,L(t)):=&\lambda_1 L(t)^{\frac{1}{2}}\sgn{w}^{\frac{1}{2}}+\mu_1M w \\
\phi_2(w,t;M,L(t)):=&\lambda_0L(t)\sign{w}\\&+\lambda_1\mu_0L(t)^{\frac{1}{2}}M\sgn{w}^{\frac{1}{2}}+\mu_0\mu_1M^2w,
\end{align}
Also consider our algorithm with $\kappa(t)$ given in Table~\mbox{\ref{tab:TVG}-(iv)} with $\alpha=1$, $\beta=0.1$ $L=1$, $M=0.1$ and $\mathcal{M}=3.2$.

The trajectory of the error, and the integral of the magnitude of the error correction term $\mathcal{H}(e_0,t;T_c)$ is shown in  the second column of Fig.~\ref{Fig:MaxErrors}. 
 Notice in the second row of Fig.~\ref{Fig:MaxErrors} that with our approach, the quality of the transient is significantly improved. Moreover, notice in the third row of Fig.~\ref{Fig:MaxErrors} that having the degree of freedom to select different classes of TVGs results in further reducing the maximum error.

 \begin{figure}
     \centering
     \def\svgwidth{8.8cm}
 \begingroup%
   \makeatletter%
   \providecommand\color[2][]{%
     \errmessage{(Inkscape) Color is used for the text in Inkscape, but the package 'color.sty' is not loaded}%
     \renewcommand\color[2][]{}%
   }%
   \providecommand\transparent[1]{%
     \errmessage{(Inkscape) Transparency is used (non-zero) for the text in Inkscape, but the package 'transparent.sty' is not loaded}%
     \renewcommand\transparent[1]{}%
   }%
   \providecommand\rotatebox[2]{#2}%
   \newcommand*\fsize{\dimexpr\f@size pt\relax}%
   \newcommand*\lineheight[1]{\fontsize{\fsize}{#1\fsize}\selectfont}%
   \ifx\svgwidth\undefined%
     \setlength{\unitlength}{1080bp}%
     \ifx\svgscale\undefined%
       \relax%
     \else%
       \setlength{\unitlength}{\unitlength * \real{\svgscale}}%
     \fi%
   \else%
     \setlength{\unitlength}{\svgwidth}%
   \fi%
   \global\let\svgwidth\undefined%
   \global\let\svgscale\undefined%
   \makeatother%
   \begin{picture}(1,1.04342712)%
     \lineheight{1}%
     \setlength\tabcolsep{0pt}%
     \put(0,0){\includegraphics[width=\unitlength,page=1]{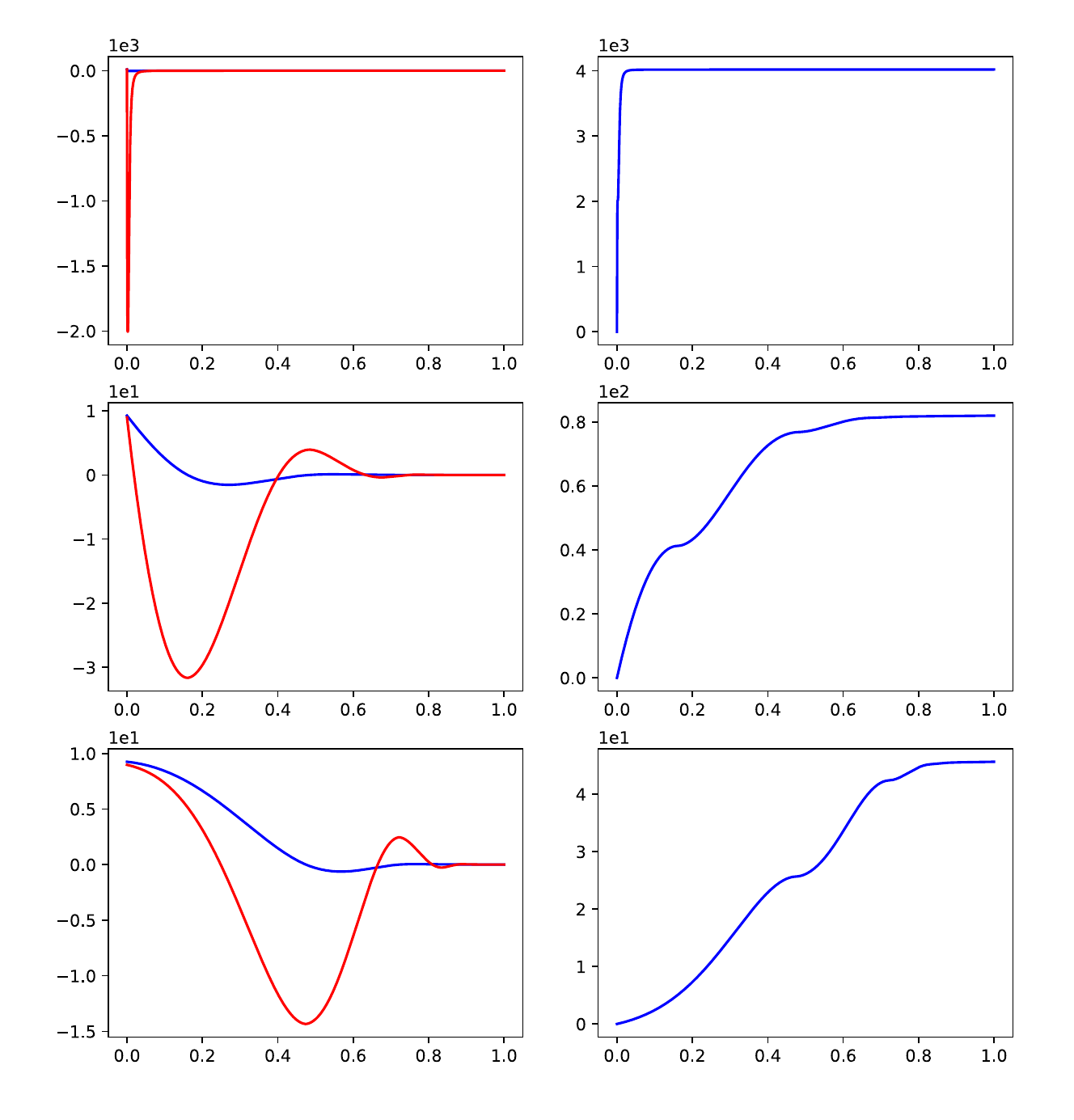}}%
     \put(0,0){\includegraphics[width=\unitlength,page=1]{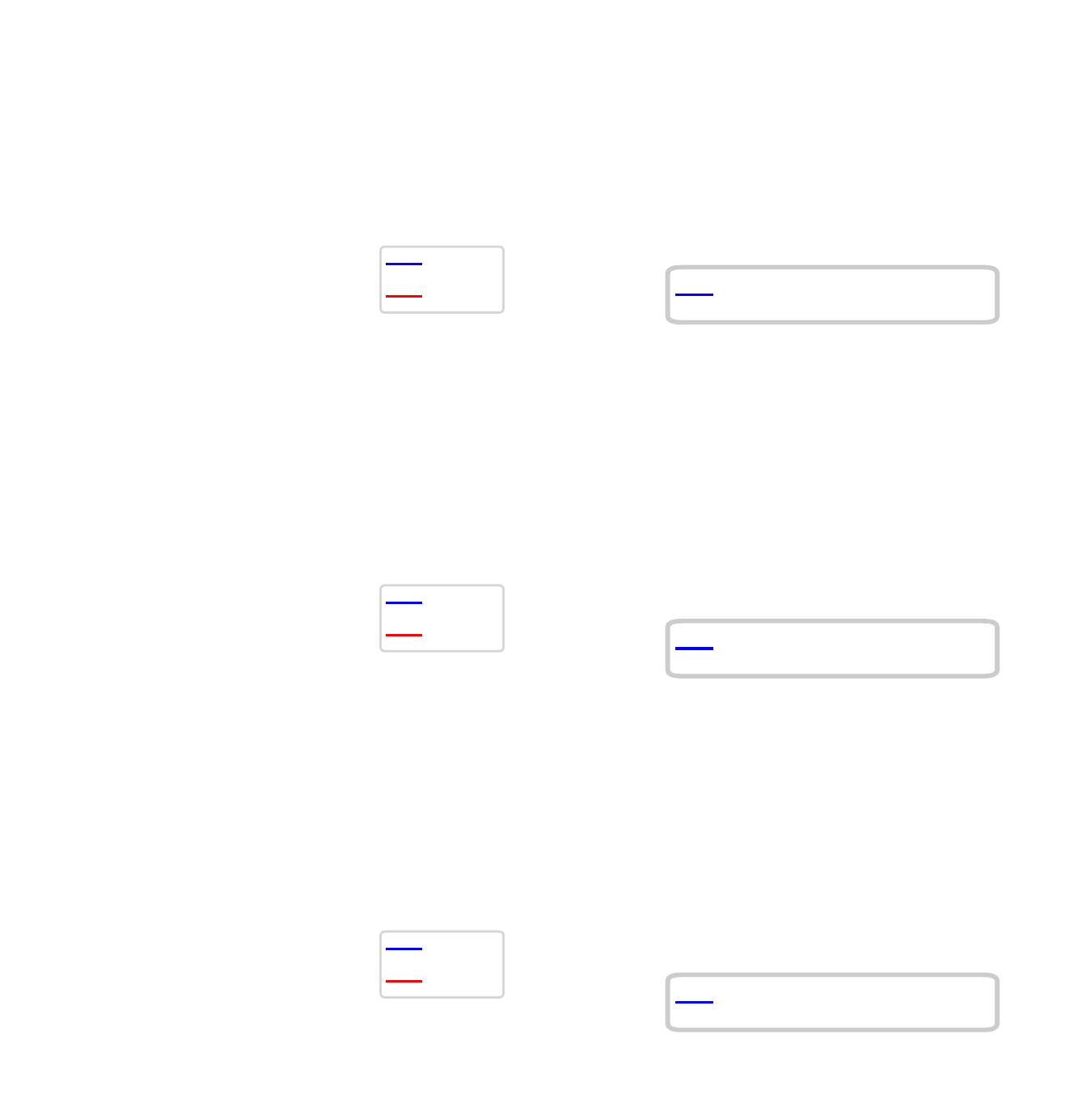}}%
     \tiny{
     \put(0.03,0.10){\color[rgb]{0,0,0}\rotatebox{90}{\makebox(0,0)[lt]{\lineheight{1.25}\smash{\begin{tabular}[t]{l}Ours- Table~\ref{tab:TVG}-(iv)\end{tabular}}}}}%
     \put(0.03,0.4){\color[rgb]{0,0,0}\rotatebox{90}{\makebox(0,0)[lt]{\lineheight{1.25}\smash{\begin{tabular}[t]{l}Ours- Table~\ref{tab:TVG}-(i)\end{tabular}}}}}%
     \put(0.03,0.72){\color[rgb]{0,0,0}\rotatebox{90}{\makebox(0,0)[lt]{\lineheight{1.25}\smash{\begin{tabular}[t]{l}\citep{Seeber2020ExactBound}\end{tabular}}}}}%
     \put(0.39649323,0.79409349){\color[rgb]{0,0,0}\makebox(0,0)[lt]{\lineheight{1.25}\smash{\begin{tabular}[t]{l}$e_0(t)$\end{tabular}}}}%
     \put(0.39649323,0.76347943){\color[rgb]{0,0,0}\makebox(0,0)[lt]{\lineheight{1.25}\smash{\begin{tabular}[t]{l}$e_1(t)$\end{tabular}}}}%
     \put(0.39649323,0.47853322){\color[rgb]{0,0,0}\makebox(0,0)[lt]{\lineheight{1.25}\smash{\begin{tabular}[t]{l}$e_0(t)$\end{tabular}}}}%
     \put(0.39649323,0.44791917){\color[rgb]{0,0,0}\makebox(0,0)[lt]{\lineheight{1.25}\smash{\begin{tabular}[t]{l}$e_1(t)$\end{tabular}}}}%
     \put(0.39649323,0.1535533){\color[rgb]{0,0,0}\makebox(0,0)[lt]{\lineheight{1.25}\smash{\begin{tabular}[t]{l}$e_0(t)$\end{tabular}}}}%
     \put(0.39649323,0.12293925){\color[rgb]{0,0,0}\makebox(0,0)[lt]{\lineheight{1.25}\smash{\begin{tabular}[t]{l}$e_1(t)$\end{tabular}}}}%
     \put(0.665,0.10423455){\color[rgb]{0,0,0}\makebox(0,0)[lt]{\lineheight{1.25}\smash{\begin{tabular}[t]{l}$\int_{0}^{t}\|\mathcal{H}(e_0,\xi;T_c)\|d\xi$\end{tabular}}}}%
     \put(0.665,0.43392438){\color[rgb]{0,0,0}\makebox(0,0)[lt]{\lineheight{1.25}\smash{\begin{tabular}[t]{l}$\int_{0}^{t}\|\mathcal{H}(e_0,\xi;T_c)\|d\xi$\end{tabular}}}}%
     \put(0.665,0.7636142){\color[rgb]{0,0,0}\makebox(0,0)[lt]{\lineheight{1.25}\smash{\begin{tabular}[t]{l}$\int_{0}^{t}\|\mathcal{H}(e_0,\xi;T_c)\|d\xi$\end{tabular}}}}%
     \put(0.26,0.03){\color[rgb]{0,0,0}\makebox(0,0)[lt]{\lineheight{1.25}\smash{\begin{tabular}[t]{l}time\end{tabular}}}}%
     \put(0.73,0.03){\color[rgb]{0,0,0}\makebox(0,0)[lt]{\lineheight{1.25}\smash{\begin{tabular}[t]{l}time\end{tabular}}}}%
     \put(0.15,1){\color[rgb]{0,0,0}\makebox(0,0)[lt]{\lineheight{1.25}\smash{\begin{tabular}[t]{l}$z_0(0)=z_1(0)=1e1$\end{tabular}}}}%
     }
   \end{picture}%
 \endgroup%
     \caption{Simulation for Example~\ref{Ex:Seeber}, with desired UBST given by $T_c=1$. In the first column: the error signals. In the second column: the integral of the magnitude of the error correction term $\int_{t_0}^{t}\|\mathcal{H}(e_0,\xi;T_c)\|d\xi$.}
     \label{Fig:MaxErrors}
 \end{figure}
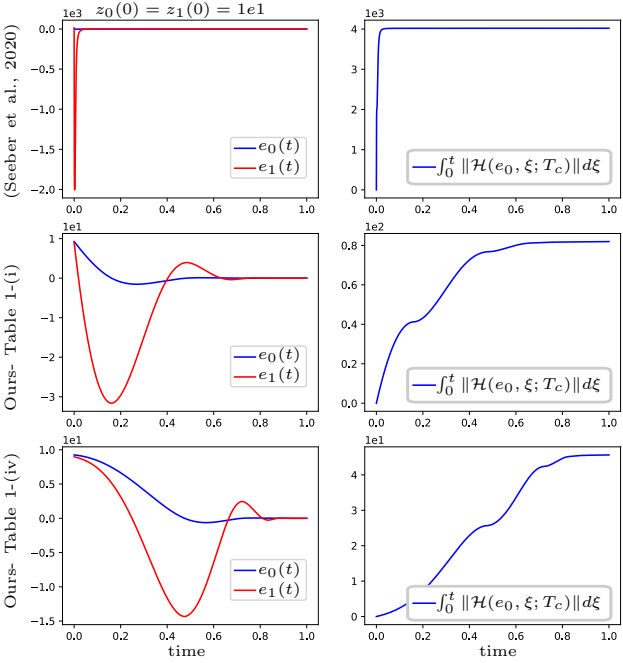

\end{example}

\section{Conclusion}
\label{Sec:Conclu}
We introduced an arbitrary-order exact differentiation algorithm with a predefined UBST for signals whose $(n+1)$-th derivative has an exponential growth bound. Compared to other differentiators based on TVGs, our approach attains a zero differentiation error before the singularity in the TVG occurs. This condition is crucial because it keeps the TVG bounded for any compact set of initial conditions. We analyze the response of our differentiator close to the predefined convergence time, finding it maintains similar robustness to measurement noise as the original differentiator by \cite{Levant2018GloballyGains}. 

In future work, we will consider the discretization of our algorithm following existing works ~\citep{Livne2014ProperDifferentiators,Carvajal-Rubio2021,Wetzlinger2019Semi-implicitDifferentiator}. Of particular interest is to obtain a consistent discretization that maintains the convergence properties of our differentiator~\citep{Koch2019Discrete-timeAlgorithm,Polyakov2019ConsistentSystems}. Additionally, we also consider extending our approach to get global ---instead of semi-global--- predefined-time convergence that is robust to measurement noise by applying the redesign methodology herein presented with a base differentiator that is fixed-time convergent~\citep{Aldana-Lopez2022AGains}.

\appendix

\section{Appendix}

\subsection{Preliminaries on time-scale transformations}

The trajectories corresponding to the system solutions of~\eqref{Eq:DiffErr} are interpreted, in the sense of differential geometry~\citep{Kuhnel2015DifferentialGeometry}, as regular parametrized curves~\citep{Pico2013}. Since we apply regular parameter transformations over the time variable, this reparametrization is referred to as time-scale transformation.

\begin{definition}(Regular parametrized curve~\cite[Definition~2.1]{Kuhnel2015DifferentialGeometry})
\label{Def:RegularParamCurve}
A regular parametrized curve, with parameter $t$, is a $C^1(\mathcal{I})$ immersion $c: \mathcal{I}\to \mathbb{R}$, defined on a real interval $\mathcal{I} \subseteq \mathbb{R}$. This means that $\frac{dc}{dt}\neq 0$ holds everywhere.
\end{definition}

\begin{definition}(Regular curve~\cite[Pg.~8]{Kuhnel2015DifferentialGeometry})
\label{Def:RegularCurve}
A regular curve is an equivalence class of regular parametrized curves, where the equivalence relation is given by regular (orientation preserving) parameter transformations $\varphi$, where $\varphi:~\mathcal{I}~\to~\mathcal{I}'$ is $C^1(\mathcal{I})$, bijective and $\frac{d\varphi}{dt}>0$. Therefore, if $c:\mathcal{I}\to\mathbb{R}$ is a regular parametrized curve and $\varphi:\mathcal{I}\to \mathcal{I}'$ is a regular parameter transformation, then $c$  and  $c\circ\varphi:\mathcal{I}'\to\mathbb{R}$ are considered to be equivalent, where $(c\circ\varphi)(t) = c(\varphi(t))$.
\end{definition}

\begin{lemma}\citep{aldana2019design}
\label{Lemma:ParTrans}
Let $\Omega(\bullet)$ satisfy Assumption~\ref{Assump:NonAut} and let $\varphi(t)$ be such $\varphi^{-1}(\tau) := T_c \int_0^\tau \Omega(\xi) d \xi$. Then, $t=\varphi^{-1}(\tau)$
is a parameter transformation (time-scaling).
\end{lemma}

\subsection{Proofs of the Main Result}
\label{Appendix:MainResult}

Before giving a proof for the main result, consider the following auxiliary result.

\begin{lemma}
\label{Lemma:F}
Let $\epsilon:=[\epsilon_0,\ldots,\epsilon_n]^T$ and let $\Omega(\bullet)$ and $\Phi(\bullet, \bullet, \bullet)$ satisfy the conditions of Theorem~\ref{Theorem:Main} and let $\delta(\tau)$ be a disturbance satisfying $\left|\delta(\tau 
    ) \right|\leq\mathcal{L}(\tau)$, and $\mathcal{M}$ such that $\exists \tau^*>0$, $\frac{1}{\mathcal{L}(\tau)}\left|\frac{d \mathcal{L}(\tau)}{d\tau}\right|\leq\mathcal{M}$ for all  $\tau>\tau^*$.
Then, the origin of the system
\begin{multline}
\label{AuxSyst}
    \frac{d\epsilon}{d\tau}=-\Phi(\epsilon_0;\mathcal{M},\mathcal{L}(\tau))+\mathcal{U}\epsilon+\mathcal{B}\delta(\tau) \\
    +\bigg(\Omega(\tau)^{-1}\frac{d\Omega(\tau)}{d\tau}+c\bigg)\mathcal{Q}(c)^{-1}\mathcal{D}\mathcal{Q}(c)\epsilon
\end{multline}
is globally finite-time stable.
\end{lemma}
\begin{pf}
Rewriting system \eqref{AuxSyst} using the coordinates $\psi(\tau) = \frac{\epsilon(\tau)}{\mathcal{L}(\tau)}$ yields
\begin{multline}
    \frac{d\psi}{d\tau} = -\Phi(\psi_0; \mathcal{M}, 1) 
    - \frac{1}{\mathcal{L}(\tau)} \frac{d\mathcal{L}(\tau)}{d\tau} \psi
    + \mathcal{U}\psi \\
    +\mathcal{B}\frac{\delta(\tau)}{\mathcal{L}(\tau)} 
    + \bigg(\Omega(\tau)^{-1}\frac{d\Omega(\tau)}{d\tau}+c\bigg)\mathcal{Q}(c)^{-1}\mathcal{D}\mathcal{Q}(c)\psi
\end{multline}
Consider first the nominal (unperturbed) part
\begin{equation}
    \frac{d\psi}{d\tau} = -\Phi(\psi_0; \mathcal{M}, 1) 
    - \frac{1}{\mathcal{L}(\tau)} \frac{d\mathcal{L}(\tau)}{d\tau} \psi
    + \mathcal{U}\psi \\
    +\mathcal{B}\frac{\delta(\tau)}{\mathcal{L}(\tau)}.
\end{equation}
Due to $|\frac{1}{\mathcal{L}(\tau)} \frac{d\mathcal{L}(\tau)}{d\tau}| \le \mathcal{M}$ and $|\frac{\delta(\tau)}{\mathcal{L}(\tau)}| \le 1$, its trajectories are solutions of the time-invariant inclusion
\begin{equation}
\label{AuxIncl}
    \frac{d\psi}{d\tau} \in -\Phi(\psi_0; \mathcal{M}, 1) + \mathcal{U}\psi + [-\mathcal{M}, \mathcal{M}] \psi + \mathcal{B} [-1, 1].
\end{equation}
According to Theorem~\ref{Th:Levant}, this inclusion is finite-time stable.
Furthermore, it follows from \citep[Lemma 2]{Levant2018GloballyGains} that it is also uniformly exponentially stable, because the constants $Q_n$, $\Delta T_n$ in that lemma do not depend on the initial time instant nor on the initial state.
Hence, according to \citep{Clarke1998AsymptoticFunctions} there exists a smooth, strong Lyapunov function $V_1 : \mathbb{R}^{n+1} \to \mathbb{R}_{+}$.
Using standard arguments, see \citep[Lemma 9.1]{Khalil2002NonlinearSystems}, convergence of all trajectories of the perturbed system \eqref{AuxSyst} to the origin is hence guaranteed, because the perturbation is linear in $\epsilon$ and its coefficient is eventually bounded by a sufficiently small constant due to Eq.~\eqref{Eq:LimRho}.

To show that also finite-time stability is maintained, consider the homogeneous approximation at the origin of the time-invariant inclusion \eqref{AuxIncl} that is obtained by setting $\mathcal{M} = 0$.
Being a special case of \eqref{AuxIncl}, this system is still finite-time stable in addition to having a homogeneity degree minus one with respect to the weights $(n+1, \ldots, 1)$.
It is straightforward to verify that the matrix $\mathcal{Q}(c)$ and, consequently, also the matrix $\mathcal{Q}(c)^{-1}\mathcal{D}\mathcal{Q}(c)$ are lower-triangular by construction.
Hence, the lowest-degree homogeneous approximation of the perturbation term $\mathcal{Q}(c)^{-1}\mathcal{D}\mathcal{Q}(c) \psi$ has a degree of at least zero.
Finite-time stability is concluded by applying Theorem 7.4 from~\cite{Bhat2005GeometricStability} and noting that this referred theorem (and its proof) stays valid in the time-varying case, as long as the lowest-degree homogeneous approximation (with degree minus one) is time-invariant and the time-varying part is uniformly bounded with respect to time.

\end{pf}

\begin{pf}[Proof of Theorem~\ref{Theorem:Main}]
Let $e_i(t)=z_i-\frac{d^iy(t)}{dt^i}$, $0,\ldots,n$ . The proof is divided into two parts. First, we will show that $e_0(t)=0$, $i=0,\ldots,n$, for $t\in [\hat{t},T_c)$ for some time $\hat{t}$. Afterwards, we show that the condition $e_i(t)=0$, $i=0,\ldots,n$ is maintained for all $t>T_c$. Let $\epsilon:=[\epsilon_0,\ldots,\epsilon_n]^T$; $e:=[e_0,\ldots,e_n]^T$, and define 
$$\mathcal{A}(c):=-\mathcal{Q}^{-1}(\mathcal{U}-\alpha \mathcal{D})^{n+1}\mathcal{B}[1,0,\dots,0].$$
Then, the dynamic for the differentiation error can be written as:
\begin{multline}
    \dot{e}=-\Lambda(t)\mathcal{Q}(c)\left[\Phi\big(e_0;\mathcal{M},\mathcal{L}(\varphi(t))\big)-\mathcal{A}(c)e\right]\\+\mathcal{U}e-\mathcal{B}y^{(n+1)}(t)
\end{multline}
With the coordinate change $\epsilon=\kappa(t)\mathcal{Q}(c)^{-1}\Lambda(t)^{-1}e,$
and considering that
$$\frac{d\kappa(t)^{-i}}{dt}\kappa(t)^{i}=-i\kappa(t)^{-1}\frac{d\kappa(t)}{dt},$$ $$\Lambda(t)^{-1}\mathcal{U}\Lambda(t)=\kappa(t)\mathcal{U} 
\text{ and } \mathcal{Q}(c)^{-1}\Lambda(t)^{-1}\mathcal{B}=\kappa(t)^{-n}\mathcal{B}.$$
Then, the dynamic of the $\epsilon$ variable, is given by 
\begin{multline}
\label{Eq:Epsilon}
    \dot{\epsilon}=\kappa(t)\bigg(\mathcal{Q}(c)^{-1}[\mathcal{U}-\kappa(t)^{-2}\frac{d\kappa(t)}{dt}\mathcal{D}]\mathcal{Q}(c)\epsilon\\-\Phi(\epsilon_0;\mathcal{M},\mathcal{L}(\varphi(t)))\\
    +\mathcal{A}(c)\epsilon-\kappa(t)^{-(n+1)}\mathcal{B}y^{(n+1)}(t)\bigg),
\end{multline}
Now, consider the time-scaling given in Lemma~\ref{Lemma:ParTrans} and notice that $\mathcal{L}(\tau):=L(\varphi^{-1}(\tau))\rho(\tau)^{-(n+1)}$, where $$\rho(\tau)=\frac{1}{T_c}\Omega(\tau)^{-1}=\left.\kappa(t)\right|_{t=\varphi^{-1}(\tau)}$$ 
and $L(\varphi^{-1}(\tau))=\left. L(t)\right|_{t=\varphi^{-1}(\tau)}$. 
Notice that, 
$$\frac{1}{\mathcal{L}(\tau)}\left| \frac{d\mathcal{L}(\tau)}{d\tau} \right|\leq \left|M\rho(\tau)^{-1}-(n+1)\rho(\tau)^{-1}\frac{d\rho(\tau)}{d\tau}\right|.$$
Thus, if $\mathcal{M}>(n+1)c$, there exists $\tau^*$ such that $\frac{1}{\mathcal{L}(\tau)}\left| \frac{d\mathcal{L}(\tau)}{d\tau} \right|\leq \mathcal{M}$, since $\rho(\tau)^{-1}$ tends to zero and $\rho(\tau)^{-1} \frac{d\rho(\tau)}{d\tau}$ tends to $c$ for $\tau \to \infty$. 

Then,
$ \frac{d\epsilon}{d\tau}=\left.  \frac{d\epsilon}{dt}  \frac{dt}{d\tau}\right|_{t=\varphi^{-1}(\tau)}.$
Since $\left.\frac{dt}{d\tau}\right|_{t=\varphi^{-1}(\tau)}= \kappa(t)^{-1}$, for $t\in[0,T_c)$, then the dynamics of~\eqref{Eq:Epsilon} in the $\tau$-time is given by
\begin{multline}
    \frac{d\epsilon}{d\tau}=\mathcal{Q}(c)^{-1}[\mathcal{U}+\bigg(\Omega(\tau)^{-1}\dfrac{d\Omega(\tau)}{d\tau}+c-c\bigg)\mathcal{D}]\mathcal{Q}(c)\epsilon\\-\Phi(\epsilon_0;\mathcal{M},\mathcal{L}(\tau))+\mathcal{A}(c)\epsilon-\mathcal{B}\delta(\tau),
\end{multline}
Since, $\mathcal{Q}(c) \mathcal{U} = (\mathcal{U} - c \mathcal{D}) \mathcal{Q}(c) + \mathcal{Q}(c) \mathcal{A}(c)$, then 
$
\mathcal{Q}(c)^{-1}[\mathcal{U}-c\mathcal{D}]\mathcal{Q}(c)=\mathcal{U}-\mathcal{A}(c),
$ and
\begin{multline}
\label{Eq:AuxSyst}
    \frac{d\epsilon}{d\tau}
    =-\Phi(\epsilon_0;\mathcal{M},\mathcal{L}(\tau))+\mathcal{U}\epsilon+\mathcal{B}\delta(\tau) \\
    +\bigg(\Omega(\tau)^{-1}\dfrac{d\Omega(\tau)}{d\tau}+c\bigg)\mathcal{Q}(c)^{-1}\mathcal{D}\mathcal{Q}(c)\epsilon,
\end{multline}
which according to Lemma~\ref{Lemma:F}, system~\eqref{Eq:AuxSyst} is finite-time stable and has a settling-time function $\mathcal{T}(\epsilon(0))$. Using Lemma~\ref{Lemma:ParTrans}, we can conclude that the settling-time function of~\eqref{Eq:DiffErr} is 
\begin{align}
\label{Eq:SeetlingTimeDif2}
    T(e(0))&=\lim_{\tau\to\mathcal{T}(\epsilon(0))}\varphi(\tau )=T_c \int_0^{\mathcal{T}(\epsilon(0))} \Omega(\xi) d \xi.
\end{align}
Thus, 
\begin{equation}
\label{Eq:STsup2}
    \sup_{e(0)\in\mathbb{R}^{n+1}}T(e(0))\leq T_c. 
\end{equation}
Then, $e_i(t)=0$, $i=1,\ldots,n$ for $t\in [\hat{t},T_c)$, where $\hat{t}=T_c\Big(1-\exp{-\alpha\mathcal{T}(\epsilon(0))}\Big)$. Moreover, it follows from~\eqref{Eq:SeetlingTimeDif2}, that the equality in~\eqref{Eq:STsup2} holds when $\sup_{e(0)\in\mathbb{R}^{n+1} }\mathcal{T}(e(0))=\infty$. Finally, note that the existence and uniqueness of Filippov solutions to \eqref{AuxSyst} for all $\tau\geq 0$ implies the existence and uniqueness of solutions for \eqref{Eq:DiffErr} for $t\in[0,T_c)$. Now, since $e(T_c)=0$ for any initial condition $e(0)$, we can continue the solution $e(t)$ for $t\geq 0$ trivially as $e(t)=0, \forall t\geq T_c$ using Theorem~\ref{Th:Levant}, because the differentiation error dynamics are given by system~\eqref{Eq:LevDiffErr}.
\end{pf}
\subsection{Proof of Proposition \ref{prop:noise}}
\label{ap:noise}
First, write the perturbed error system as $\dot{e}=-\mathcal{H}(e_0+\eta_0,t;T_c)+\mathcal{U}e-\mathcal{B}y^{(n+1)}(t)$ which is equivalent to 
\begin{multline}
\label{eq:AuxSystPert}
    \frac{d\epsilon}{d\tau}=-\Phi(\epsilon_0+\eta_0;\mathcal{M},\mathcal{L}(\tau))+\mathcal{U}\epsilon+\mathcal{B}\delta(\tau) \\
    +\bigg(\Omega(\tau)^{-1}\frac{d\Omega(\tau)}{d\tau}+c\bigg)\mathcal{Q}(c)^{-1}\mathcal{D}\mathcal{Q}(c)\epsilon
\end{multline}
under the change of coordinates $\epsilon=\kappa(t)\mathcal{Q}(c)^{-1}\Lambda(t)^{-1}e$ and $\tau=\varphi(t)$. Hence, the noise bound is written as $|\eta_0|\leq \eta \rho(T_f)^{n+1} \mathcal{L}(\varphi(\tau))$ where $T_f=\varphi(T_c^*)$. By similar arguments as in the proof of Lemma \ref{Lemma:F}, the terminal bounds in \cite{Levant2018GloballyGains} for the perturbed version of \eqref{Eq:LevDiff} apply to \eqref{eq:AuxSystPert} as $|\epsilon_i(\tau)|\leq \gamma_i{L}(\varphi^{-1}(\tau)){\rho(\tau)}^{-(n+1)}\rho(T_f)^{{n-i+1}} \eta^{\frac{n-i+1}{n+1}}$ for $\tau, T_f>0$ sufficiently large, $\tau\leq T_f$ and sufficiently small $\eta>0$. Now, in the $e$ coordinates
$$
\begin{aligned}
&|e_i| = \left|\kappa(t)^{i}\sum_{j=0}^nq_{ij}(c)\epsilon_j\right|\\
&\leq \sum_{j=0}^n \gamma_j{L}(t)q_{ij}(c)\frac{\rho(T_f)^{{n-j+1}}}{\kappa(t)^{n-i+1}} \eta^{\frac{n-j+1}{n+1}}\\
&\leq\tilde{\gamma}_iL(t)\frac{\rho(T_f)^{n+1}}{\rho(0)^{n-i+1}} \eta^{\frac{n-i+1}{n+1}}
\end{aligned}
$$
where $q_{ij}(c)$ are the components of $\mathcal{Q}(c)$ with $q_{ij}(c)=0, j>i$, and $\tilde{\gamma}_i>0$ is a constant for which it is complied that $\sum_{j=0}^n \gamma_jq_{ij}(c)\eta^{\frac{n-j+1}{n+1}}\leq \tilde{\gamma}_i\eta^{\frac{n-i+1}{n+1}}$, which always exist for sufficiently small $\eta$.

\end{document}